\documentclass[12pt]{amsart}

\usepackage{graphicx}
\usepackage{float}
\usepackage{amssymb}
\usepackage{epstopdf}
\usepackage{amsmath,amsthm,amscd,amssymb}
\usepackage{latexsym}
\usepackage[colorlinks,citecolor=red,pagebackref,hypertexnames=false]{hyperref}
\usepackage{geometry}                % See geometry.pdf to learn the layout options. There are lots.
\usepackage[us, 12hr]{datetime}

\numberwithin{equation}{section}

\theoremstyle{plain}
\newtheorem{theorem}{Theorem}[section]
\newtheorem{lemma}[theorem]{Lemma}

\theoremstyle{definition}
\newtheorem{definition}[theorem]{Definition}

\newtheorem{case[theorem]}{Case}

\theoremstyle{remark}
\newtheorem{remark}[theorem]{Remark}

\numberwithin{equation}{section}

\begin{document}
\title{
APPROXIMATE ORTHOGONALITY, BOURGAIN'S PINNED DISTANCE THEOREM AND EXPONENTIAL FRAMES} 
 %    Information for first author

\author{Alex Iosevich and Azita Mayeli}

\date{}
%\date{Last edit: \currenttime, \today}
\address{Department of Mathematics, University of Rochester, Rochester, NY}
\email{iosevich@math.rochester.edu}
\address{City University of New York, NY} 
\email{amayeli@gc.cuny.edu}
%\address{Department of Mathematics, University of Rochester, Rochester, NY}
%\email{jonpak@math.rochester.edu}

\thanks{The research of the first listed author was supported in part by the National Science Foundation under grant no. HDR TRIPODS - 1934962 and by the NSF DMS - 2154232}

\begin{abstract} Let $A$ be a countable and discrete subset of ${\Bbb R}^d$, $d \ge 2$, of positive upper Beurling density. Let $K$ denote a bounded symmetric convex set with a smooth boundary and everywhere non-vanishing Gaussian curvature. It is known that  ${\mathcal E}(A)=\{e^{2 \pi i x \cdot a}\}_{a \in A}$ cannot serve as an orthogonal basis for $L^2(K)$ \cite{IKT01}. In this paper, we prove that even approximate average orthogonality is an obstacle to the existence of an exponential frame in the following sense. Let   $A$ be  as above  and $\phi \ge 0$ be a continuous monotonically non-increasing function on $[0, \infty)$ such that the {\it approximate orthogonality condition} 
holds 
\begin{equation}\notag
{\left( \frac{1}{2^j} \int_{2^j}^{2^{j+1}} \phi^p(t) dt \right)}^{\frac{1}{p}} \leq c_j 2^{-j\frac{d+1}{2}} \quad \text{and} \quad \ |\widehat{\chi}_K(a-a')| \leq \phi(\rho^*(a-a'))\ \  \forall  a \not=a , a,a' \in A, 
\end{equation} 
 where $\rho^{*}$ is the Minkowski functional on $K^{*}$, the dual body of $K$. Then, if   
$$\limsup_{j \to \infty} c_j=0,$$ then the upper density of $A$ is equal to $0$, hence ${\mathcal E}(A)={\{e^{2 \pi i x \cdot a} \}}_{a \in A}$ {\bf is not} a frame for $L^2(K)$. The case $p=\infty$ was previously established by the authors of this paper in \cite{IM2020}. The point is that if ${\mathcal E}(A)$ is a frame for $L^2(K)$, then very few pairs of distinct exponentials $e^{2 \pi i x \cdot a}$, $e^{2 \pi i x \cdot a'}$ from ${\mathcal E}(A)$ come anywhere near being orthogonal. Our proof uses a generalization of Bourgain's result on {\bf pinned} distances determined by sets of positive Lebesgue upper density in ${\Bbb R}^d$, $d \ge 2$. We also improve the $L^{\infty}$ version of this result originally established in \cite{IM2020}. By using an extension of the combinatorial idea from \cite{IR03}, we prove that under the $L^{\infty}$ hypothesis, $A$ is finite if $d \not=1 \mod 4$. If $d=1$ mod $4$, $A$ may be infinite, but if it is, then it must be a subset of a line. 
\end{abstract} 

\maketitle

\section{Introduction} 
The primary goal  of this paper is to continue the investigation of the structure of discrete sets $A$ such that ${\{e^{2 \pi i x \cdot a}\}}_{a \in A}$ is a frame for $L^2(K)$, where $K$ is a bounded symmetric convex body with a smooth boundary and non-vanishing curvature. Recall that such a collection of exponentials is a frame if there exist $0<C_L \leq C_U<\infty$ so that for every $f \in L^2(K)$, 
\begin{equation}\label{frame} C_L|K| \cdot {||f||}^2_{L^2(K)} \leq \sum_{a \in A} {\left|\widehat{f}(a)\right|}^2 \leq C_U|K| \cdot {||f||}^2_{L^2(K)}.\end{equation}

Let us recall some definitions. 

\begin{definition} Let $S$ be a separated subset of ${\Bbb R}^d$, in the sense that there exists $\delta>0$ such that $|s-s'| \ge \delta>0$ for all $s,s' \in S$, $s \not=s'$. The upper density of $S$, is given by 
$$ \overline{dens}(S)=\limsup_{R \to \infty} \frac{\# \{S \cap {[-R,R]}^d\}}{2^dR^d}.$$ 

Similarly, the lower density of $S$ is given by 
$$ \underline{dens}(S)=\liminf_{R \to \infty} \frac{\# \{S \cap {[-R,R]}^d\}}{2^dR^d}.$$

If $\overline{dens}(S)=\underline{dens}(S)$, then we say that the Beurling density $dens(S)$ exists and is equal to this quantity. 
\end{definition} 

It is not difficult to construct a frame for $L^2(K)$ by enclosing $K$ in a rectangular box and restricting the orthonormal exponential basis for this box to $K$.  In this case, the frame is  tight, i.e.,  in the definition of the frame in \eqref{frame} we have $C_L=C_U$.  One unfortunate feature of such a frame is that its density is much higher than the volume of $K$.

\vskip.125in 

Also recall that ${\mathcal E}(A)={\{e^{2 \pi i x \cdot a}\}}_{a \in A}$ is a Riesz basis for $L^2(K)$ if it is equivalent (via an invertible and bounded map) to an orthonormal basis (not necessarily exponential) for $L^2(K)$. It is an old open question of whether $L^2(B_d)$, $d \ge 2$, possesses a Riesz basis of exponentials, where $B_d$ denotes the unit ball.

It is well-known that if ${\mathcal E}(A)$ is a Riesz basis for $L^2(K)$, then $dens(A)=|K|$. In hopes of shedding some light on the question of a possible Riesz basis for $L^2(K)$, where $K$ is bounded, symmetric, convex and has a smooth boundary with non-vanishing curvature, we prove a result that puts a significant restriction on the structure of a set $A$ that generates a frame for $L^2(K)$. 

A well-known stationary phase formula due to Herz \cite{H62} (see also \cite{GS58} and Lemma \ref{stationaryphase} below for a precise formulation) says that 

\begin{equation} \label{firstterm} \widehat{\chi}_K(\xi)=\kappa^{-\frac{1}{2}} \left( \frac{\xi}{|\xi|} \right) 
\sin \left( 2 \pi \left( \rho^{*}(\xi)-\frac{d-1}{8}\right) \right){|\xi|}^{-\frac{d+1}{2}} \end{equation} plus an error term,  where $\rho^{*}$ is the Minkowski functional on $K^{*}$, the dual body of $K$. We are going to show in Theorem \ref{maintool} below that if ${\mathcal E}(A)$ is a frame, $|\widehat{\chi}_K(a-a')|$, with $a \not=a'$, $a,a' \in A$,  cannot be much smaller than ${|a-a'|}^{-\frac{d+1}{2}}$ very often. In other words, the difference set $A-A$ must stay away from the zeroes of $\widehat{\chi}_K$. 

%\footnote{If the exponentials are a frame, then by the assumption the set A is well-sepearated. Assume that $\delta$ is the infimum of the separation distance. In this case, $|\widehat{\chi}_K(a-a')|$ can not be less than $\delta^{-\frac{d+1}{2}}$ more often. This means that for most part of $A-A$, the values $|\widehat{\chi}_K(a-a')|$ are bounded from below by $\delta$. I do not know if this view helps somewhere here.}

It was shown in \cite{IKT01} (see also the case of the ball in \cite{IKP01}) that with $K$ as above, $L^2(K)$ does not possess an orthogonal basis of exponentials. Orthogonality implies that $(A-A) \backslash \{ \vec{0}\}$ is strictly contained in the zero set of $\widehat{\chi}_K$. What we prove in this paper is that $A-A$ cannot be close to the zeroes of $\widehat{\chi}_K$ very often, even on average, in order for ${\mathcal E}(A)$ to be a frame for $L^2(K)$. More precisely, the main  results of the current paper   are  the following. 

\subsubsection{Our Contributions} 
\begin{theorem} \label{main} Let $A$ be a countable and discrete subset of ${\Bbb R}^d$, $d \ge 2$, and let $ p \in [1, \infty]$. Let $K$ denote a bounded symmetric convex set with a smooth boundary and everywhere non-vanishing Gaussian curvature. Let $\phi \ge 0$ be a continuous monotonically non-increasing function on $[0, \infty)$ such that 
\begin{equation} \label{boundfunctionprop} {\left( \frac{1}{2^j} \int_{2^j}^{2^{j+1}} \phi^p(t) dt \right)}^{\frac{1}{p}} \leq c_j \left(2^{-j\frac{d+1}{2}}\right) \ \  \ \forall j\in \Bbb Z,  \end{equation}
for some constants $c_j\geq 0$. 

Further suppose that 
\begin{equation} \label{boundfunction} |\widehat{\chi}_K(a-a')| \leq \phi(\rho^{*}(a-a')) \end{equation}  for all $a \not=a'$, $a,a' \in A$. 

\vskip.125in 

Then if 
$$\limsup_{j \to \infty} c_j=0,$$ the upper density of $A$ is equal to $0$, and, consequently, ${\mathcal E}(A)$ {\bf is not} a frame for $L^2(K)$. 
\end{theorem} 

\vskip.125in 

\begin{remark} Note that if ${\mathcal E}(A)$ is an orthogonal system, then $\widehat{\chi}_K(a-a')  = 0$ if $a \not=a'$. Therefore, Theorem \ref{main} is a generalization of the fact   that $L^2(K)$ does not possess any orthogonal basis of exponentials if $K$ is a symmetric convex body with a smooth boundary and non-vanishing Gaussian curvature \cite{IKT01}.  \end{remark} 

\begin{remark} An immediate consequence of Theorem \ref{main} is that if ${\mathcal E}(A)$ is a frame for $L^2(K)$, with $K$ as above, then for 
any subset $A'\subset A$ on  which (\ref{boundfunctionprop}) and (\ref{boundfunction}) hold, with $\limsup_{j \to \infty} c_j=0$,    $A'$ has zero  positive density. In other words, these properties can hold only on a subset of zero density. It would be interesting to determine if a stronger statement is true, namely that if  ${\mathcal E}(A)$ is a frame for $L^2(K)$, and the above hypotheses hold on $A' \subset A$, then for some $0<\alpha<d$, 
$$ \limsup_{R \to \infty} \frac{\# \{A' \cap B_R \}}{R^{d-\alpha}}=0$$ 
where $B_R$ denote the ball of radius $R$.  
 For example, if $A$ were contained in a line, we could take any $\alpha<d-1$.  We partially address this question in Theorem \ref{inftyforrealth}. \end{remark} 

\vskip.125in 

\begin{remark} It is proved in (\cite{IP00}; see also \cite{IK06}) that if $\Omega$ is a bounded domain of positive Lebesgue measure in ${\Bbb R}^d$, $d \ge 2$, and ${\{e^{2 \pi i x \cdot a}\}}_{a \in A}$ is a frame for $L^2(\Omega)$ with upper and lower frame constants $C_U, C_L$, respectively, then every  $d$-dimensional cube of side-length $\ell$ with 
$$ \ell>{\left( \frac{C_U {|\partial \Omega|}_{\alpha}}{C_L |\Omega|} \right)}^{\frac{1}{d-\alpha}}$$ contains at least one point of $A$, where $0<\alpha<d$ and ${|\partial \Omega|}_{\alpha}$ denotes the $\alpha$-dimensional Minkowski content of $\partial \Omega$, the boundary of $\Omega$. In particular, if $\Omega$ is convex, as it is assumed in Theorem \ref{main}, every cube of side-length  $\ell$ with 
$$ \ell > \frac{C_U |\partial \Omega|}{C_L |\Omega|} $$

%\footnote{I wonder what becomes of this inequality when the frame is orthogonal. Is it also sufficient? Can the same results be proved using the sinc function on the time domain? I did not realize it before how nice this result is.} 

contains at least one point of $A$, where $|\partial \Omega|$ denotes the surface measure of $\partial \Omega$. This can be viewed as a quantitative version of the fact (due to Beurling (\cite{B66})), which states  that if $\Omega$ has positive Lebesgue measure and ${\mathcal E}(A)$ is frame for $L^2(\Omega)$, then $A$ has positive upper density. \end{remark} 

\vskip.125in 

\begin{remark}\label{rem} The case $p=\infty$, $K=B_d$, the unit ball in ${\Bbb R}^d$, is essentially contained in the prequel to this paper (\cite{IM2020}). More precisely, the authors showed that if $\phi: [0, \infty)\to [0, \infty)$ is any bounded and positive measurable  function such that 
\begin{equation} \label{previous} \lim_{t \to \infty} {(1+t)}^{\frac{d+1}{2}} \phi(t)=0,    \ \ \text{and} \quad     |\widehat{\chi}_{B_d}(a-a')| \leq \phi(|a-a'|)  \quad   \forall   \ a,a' \in A, \  a\neq a', \end{equation}
 then $A$ has zero upper density. Making a change of variables $s=\frac{t}{2^j}$ in (\ref{boundfunctionprop}),  it turns the inequality  (\ref{boundfunctionprop})  into 
$$ {\left( \int_1^2 \phi^p(2^js) ds \right)}^{\frac{1}{p}} \leq c_j 2^{-j\frac{d+1}{2}}.$$ 

The $L^{\infty}$ version of this statement is 
$$ \sup_{s \in [1,2]} 2^{j\frac{d+1}{2}} \phi(2^j s) \leq c_j,$$ so in the case $p=\infty$, (\ref{boundfunctionprop}) is a slightly weaker assumption than the (\ref{previous}), except that the monotonicity assumption is not needed in the $L^{\infty}$ setting. However, as the reader is about to see, we can strengthen the $L^{\infty}$ result considerably by using different methods. 
\end{remark} 

\vskip.125in 

Our second result is a significant improvement of a previous theorem by the authors of this paper described  in Remark \ref{rem} above. 

\begin{theorem} \label{inftyforrealth} Suppose that if $\phi: [0, \infty)\to [0, \infty)$ is any bounded positive measurable function such that (\ref{previous}) holds for all $a,a' \in A$, $a \not=a'$. Then $A$ is contained in a line. Moreover, if $d \not=1 \mod 4$, then $A$ is finite. \end{theorem} 

\vskip.125in 

\begin{remark} As the reader will see, the proof of Theorem \ref{inftyforrealth} is based on the techniques developed in \cite{IR03}. The re-examantion of these techniques in the context of Theorem \ref{inftyforrealth} allowed the authors of this paper to realize that the approximate version of the Erd\H os Integer Distance Principle proved in \cite{IR03} holds under considerably weaker hypotheses. \end{remark} 

\vskip.25in 

\section{Outline of the proofs and the structure of the paper} 

\vskip.125in 

In this section we outline the proof of Theorems \ref{main} and  \ref{inftyforrealth}. 

\vskip.125in 

\subsection{Outline of the proof of Theorem \ref{main}} By the method of stationary phase (see Lemma \ref{stationaryphase} below), with $K$ as above, 
$$|\widehat{\chi}_K(a-a')| \leq C{(1+|a-a'|)}^{-\frac{d+1}{2}}$$ for some constant $C>0$. It follows that (\ref{keyaverage}) always holds with $c=C$. The point is that if it holds with a small enough constant, we can invoke the asymptotic formula (\ref{herzformula}) (see below) and argue that 
$$ \sin \left(2 \pi \left(\rho^{*}(a-a')+\frac{d-1}{8}\right)\right)$$ is small quite often for $a \not=a'$, $a,a' \in A$. This will be accomplished by reducing the proof of Theorem \ref{main} to the following statement. 

\begin{theorem} \label{maintool} Let $A$ be a discrete and countable subset of ${\Bbb R}^d$, $d \ge 2$. Let $K$ denote a bounded symmetric convex set with a smooth boundary and everywhere non-vanishing Gaussian curvature. Let ${\mathcal A}_j$ denote annuli of inner radius $2^j$, and outer radius $2^{j+1}$, centered at $a' \in A$. Suppose that 
\begin{equation} \label{keyaverage} 
{\left( \frac{1}{2^{dj}} \sum_{a \in A \cap {\mathcal A}_j} {|\widehat{\chi}_{K}(a-a')|}^p \right)}^{\frac{1}{p}} \leq c_j 2^{-j\frac{d+1}{2}} \end{equation} with $c_j$ independent of $a'$. Then if 
\begin{equation}\label{limsuptozero}
\limsup_{j \to \infty} c_j=0,\end{equation}
then the upper density of $A$ is equal to $0$, hence ${\mathcal E}(A)$ {\bf is not} a frame for $L^2(K)$. 
\end{theorem} 

%\footnote{all the $a$ in the intersection have almost the same distance from $a'$. If we assume that $|a-a'|= o(1)$, and use the inequality for $|\hat \chi_K(a-a')|$, then the inequality (\ref{keyaverage}) can be estimated from above using the $\sharp A\cap \mathcal A_j$. Comparing this with the right hand side of (2.1) we may get a sharper inequality for large j.}

\vskip.125in 

By pigeonholing, we will be able to conclude that there exists a point $x_0$ living in a small neighborhood of an $a_0 \in A$, such that pinned $\rho^{*}$-distances from $x_0$ to the small neighborhood of a positive proportion of $A$ is clustered around $\frac{k}{2}+\frac{d-1}{8}$, where $k$ ranges over the integers. This violates a variant of the result due to Bourgain which says that in sets of positive upper Lebesgue density, a pinned distance set with respect to almost every point contains every sufficiently large positive real number. In order to invoke Bourgain's result we will further have to refine the previously pigeon-holed subset of $A$. 

\vskip.125in 

\subsection{Outline of the proof of Theorem \ref{inftyforrealth}} 

\vskip.125in 

The authors of  \cite{IR03}  proved that if ${\mathcal E}(A)$ is an infinite orthogonal set in $L^2(K)$, where $K\subset\Bbb R^d$ is a symmetric convex body with a smooth boundary and everywhere non-vanishing curvature, then $A$ is contained in a line, thus, ${\mathcal E}(A)$ is incomplete in $L^2(K)$.   Moreover, they showed that if $d \not=1 \mod 4$, then $A$ must be finite. We are going to see that the proof of this result still goes through under the weaker assumption given by (\ref{previous}) above. 

\vskip.125in 

\subsection{Structure of this paper} This paper is organized as follows. In the Section \ref{proofofmaintool} below, we prove Theorem \ref{maintool}. At the end of that section, we comment on the modification of the pinned distance set result due to Bourgain that is needed for our proof. In Section \ref{mainfrommaintool} we deduce Theorem \ref{main} from Theorem \ref{maintool}. 

\vskip.125in 

\section{Proof of Theorem \ref{maintool}} 
\label{proofofmaintool} 

\vskip.125in 

Suppose that (\ref{keyaverage}) and (\ref{limsuptozero}) hold and $A$ has finite upper density $\alpha_0>0$. By assumption, for a fixed $a'\in A$ 
\begin{equation} \label{pigeon} \sum_{a \in A \cap {\mathcal A}_j} {|\widehat{\chi}_{K}(a-a')|}^p \leq c^p_j2^{j(d-p\frac{d+1}{2})}. \end{equation} We shall need the following well-known stationary phase result. 

\begin{lemma}\label{stationaryphase} (\cite{H62}) Let $K$ be a bounded symmetric convex body with a smooth boundary and everywhere non-vanishing Gaussian curvature. Let 
$$ \rho^{*}(\xi)=\sup_{x \in K} x \cdot \xi.$$ 

Given $\omega \in S^{d-1}$, let $\kappa(\omega)$ denote the Gaussian curvature of $\partial K$ at the (unique) point where the unit normal is $\omega$. Then 
\begin{equation} \label{herzformula} \widehat{\chi}_K(\xi)=\kappa^{-\frac{1}{2}} \left( \frac{\xi}{|\xi|} \right) 
\sin \left( 2 \pi \left( \rho^{*}(\xi)-\frac{d-1}{8}\right) \right){|\xi|}^{-\frac{d+1}{2}}+{\mathcal D}_K(\xi), \end{equation} where 
$$ |{\mathcal D}_K(\xi)| \leq C_K {|\xi|}^{-\frac{d+3}{2}}. $$ 

%and 
%$|\kappa^{-\frac{1}{2}} \left( \frac{\xi}{|\xi|} \right) 
%| =c$ \footnote{is this true?} 
\end{lemma} 

Plugging this into (\ref{pigeon}) yields 

\begin{equation} \label{errorabsorbed} \sum_{a \in A \cap {\mathcal A}_j} \kappa^{-\frac{p}{2}} \left( \frac{a-a'}{|a-a'|} \right) 
{\left|\sin \left( 2 \pi \left( \rho^{*}(a-a')-\frac{d-1}{8}\right) \right)\right|}^p {|a-a'|}^{-p\frac{d+1}{2}} \leq c^p_jC_12^{j(d-p\frac{d+1}{2})} \end{equation} for some fixed $C_1>0$, and it follows that 

$$ \sum_{a \in A \cap {\mathcal A}_j}  {\left| \sin \left( 2 \pi \left( \rho^{*}(a-a')-\frac{d-1}{8}\right) \right)\right|}^p 
\leq C_2 c^p_j 2^{dj},$$ for some fixed $C_2>0$ that depends on  the maximum value of the curvature function on $\partial K$. Note that the error term 
${\mathcal D}_K$ from (\ref{herzformula}), which is $O(2^{-j \frac{d+3}{2}})$ when $a \in {\mathcal A}_j$,  can be absorbed into the constant $C_1$ on the right hand side of (\ref{errorabsorbed}) if $j$ is sufficiently large. 

\vskip.125in 

We write the left hand side in the form $I+II$, where $I$ is the sum over 
$$ \left\{a \in A \cap {\mathcal A}_j: 
{\left| \sin \left( 2 \pi \left( \rho^{*}(a-a')-\frac{d-1}{8}\right) \right)\right|}^p \leq \delta \right\}$$ and $II$ is the sum over the complement of this set, with $0<\delta<<1$ to be determined later. We have 

$$ II> \delta \cdot \# \left\{a \in A \cap {\mathcal A}_j: 
{\left|\sin \left( 2 \pi \left( \rho^{*}(a-a')-\frac{d-1}{8}\right) \right)\right|}^p>\delta \right\},$$ hence 

\begin{align}\label{estimate-x}
\# \left\{a \in A \cap {\mathcal A}_j: 
{\left|\sin \left( 2 \pi \left( \rho^{*}(a-a')-\frac{d-1}{8}\right) \right)\right|}^p>\delta \right\}< C_3 c^p_j \delta^{-1} 2^{dj},
\end{align}
 for some fixed $C_3>0$. 

Suppose that $j \ge j_0$, to be determined later. If $\sup_{j \ge j_0} c^p_j$ is sufficiently small, choosing $\delta=\frac{1}{1000}$ we have that the right hand side of (\ref{estimate-x}) is $<(1-r) \# \{A \cap {\mathcal A}_j\}$, where $r \in (0,1)$, close to $1$, to be chosen later. This can be done since  $A$ has positive upper density and 
$$\limsup_{j \to \infty} \frac{\# \{A \cap {\mathcal A}_j \} }{|{\mathcal A}_1|2^{dj}}=\alpha_0>0.$$ 

Note that the upper density is not affected by taking only the annuli ${\mathcal A}_j$ with $j \ge j_0$ for any finite $j_0$. More precisely, the upper density of $A$ with $\cup_{j \leq j_0} \{ {\mathcal A}_j \cap A\}$ removed is the same as the upper density of $A$, 

\vskip.125in 

It follows that if we throw these points of $A$ which  belong to the set in (\ref{estimate-x}) out, we will have an $r$ fraction of the points of 
$A \cap {\mathcal A}_j$ left. Doing this for every $j \ge j_0$, with $j_0$ to be determined, and taking the union, we obtain a set $A_{good,j_0}(a')$ of positive upper density, such that 
\begin{equation} \label{sol} {\left| \sin \left( 2 \pi \left( \rho^{*}(a-a')-\frac{d-1}{8}\right) \right)\right|}^p \leq \delta \ \text{for all} \ a \in A_{good, j_0}(a'). \end{equation}

It follows from (\ref{sol}) that if $E$ denote the $\delta$-neighborhood of $A_{good, j_0}(a')$, with $\delta$ as above, then there does not exist $L_0$ such that for every $L>L_0$, $\rho^{*}(x-a')=L$. If $a'$ is replaced by any point in the ball of radius $\delta$ centered at $a'$, the conclusion remains the same. We shall use a generalization of the following result due to Bourgain. 

\begin{theorem} \label{bourgain} Let $E$ be a set of positive upper Lebesgue density in ${\Bbb R}^2$. Then for almost every $x_0 \in E$, there exists $L_0=L_0(x_0)>0$ such that for all $L>L_0$, there exists $x \in E$ with $|x-x_0|=L$. 
%\footnote{ \color{red} On the side, it would be interesting to obtain this result using NN. We fix a set $E$. The inputs will be $x_0$ and the parameters defining the NN. The output is $L_0$. I will write more on details.   }
\end{theorem}  

As we shall sketch below in Section \ref{bourgainth}, Bourgain result can be generalized as follows. Let $E$ be a set of positive upper Lebesgue density in ${\Bbb R}^d$, $d \ge 2$. Let $K$ and $\rho^{*}$ be as above. Then for almost every $x_0 \in E$, there exists $L_0>0$ such that for all $L>L_0$, there exists $x \in E$ with $\rho^{*}(x-x_0)=L$. 

This result was established by Bourgain in \cite{Bou86} (see a related result in \cite{FKW90}). However, as we shall indicate below, the generalization follows from Bourgain's method with only minor modifications. The proof only relies on the $L^p$ estimate for the circular maximal function, which Bourgain proved in 1986 (\cite{Bou86II}). However, the same estimate holds if the circle is replaced with any smooth curve with non-vanishing curvature. Since  
$$ \{x \in {\Bbb R}^d: \rho^{*}(x)=1 \}$$ is a smooth surface, with $\rho^{*}$ as above, the two-dimensional argument follows. In higher dimensions, the argument is much easier as the needed maximal estimate was proved by Stein (\cite{St73}) in the case of the sphere, and by Greenleaf (\cite{G81}) for general compact surfaces with non-vanishing Gaussian curvature. We shall sketch the proof of this in Section \ref{bourgainth} below. 

 \vskip.125in 
 
In order to use Bourgain's result, we need to make sure that if $E$ is as in the previous paragraph, then one of the points in the $\frac{1}{1000}$ neighborhood of $a'$ can serve as a pin in Theorem \ref{bourgain}. The proof of Bourgain's theorem yields a slightly stronger results than the conclusion of Theorem \ref{bourgain}, namely that there exists $x_0 \in E$ such that for every subset of $E$ of relative density $r$ sufficiently close to $1$, there exists $L_0>0$ such that for all $L>L_0$, there exists $x \in E$ with $\rho^{*}(x-x_0)=L$. In view of this, we complete the proof by finding such an $x_0$ in the $\delta$-neighborhood of $A$ and choosing $a'$ to be the point in $A$ in the $\delta$-neighborhood of which $x_0$ resides. 

\vskip.25in 
 
\section{Proof of  Theorem \ref{main}} 
\label{mainfrommaintool} 

\vskip.125in 

Let ${\mathcal A}_j$ be as above. By the assumption of Theorem \ref{main},  we have 
$$\sum_{a \in  A\cap {\mathcal A}_j} {|\widehat{\chi}_{K}(a-a')|}^p \leq \sum_{a \in A \cap {\mathcal A}_j} \phi^p(\rho^{*}(a-a')).$$

Since $\phi$ is monotonic, the right expression above is bounded by  
$$ C \int_{c_12^j \leq \rho^{*}(u-a') \leq c_22^{j+1}} \phi^p(\rho^{*}(u-a')) du=C \int_{c_12^j \leq \rho^{*}(u) \leq c_22^{j+1}} \phi^p(\rho^{*}(u)) du,$$
or 
$$  C \int_{c_12^j \leq \rho^{*}(r \omega) \leq c_22^{j+1}} \phi^p(\rho^{*}(r \omega)) r^{d-1} d\omega dr=  C \int_{c_12^j \leq r \rho^{*}(\omega) \leq c_22^{j+1}} \phi^p(r\rho^{*}(\omega)) r^{d-1} d\omega dr,$$
where $C$ is an independent constant. 
\vskip.125in 

Let $t=r \rho^{*}(\omega)$. We obtain 
$$ C \int_{c_12^j \leq t \leq c_22^{j+1}} \int_{S^{d-1}} \phi^p(t) t^{d-1} \frac{d\omega}{{(\rho^{*}(\omega))}^d} dt $$
\begin{equation} \label{almost}=C \int_{c_12^j \leq t \leq c_22^{j+1}} \int_S \phi^p(t) t^{d-1} dt \frac{d\sigma_{\rho^{*}}(\eta)}{|\nabla \rho^{*}(\eta)|}, \end{equation} where $d\sigma_{\rho^{*}}$ is the surface measure on 
$$ S=\{x: \rho^{*}(x)=1\},$$ the $\rho^{*}$-unit sphere. 

The expression (\ref{almost}) follows from the classical co-area formula. See, for example, \cite{IS96}, for a similar calculation. 

Since 
$\rho^{*}(\eta)=\eta \cdot \nabla \rho^{*}(\eta),$ it follows that $|\nabla \rho^{*}(\eta)| \ge C>0$, when $\eta$ is on the $\rho^{*}$-unit sphere, for some different $C>0$. It follows that (\ref{almost}) is bounded by 
\begin{align}\notag
   (\ref{almost})  \leq C' 2^{j(d+1)} \cdot 2^{-j} \int_{c_12^j}^{c_2 2^{j+1}} \phi^p(t) t^{d-1} dt 
  \leq C'' c_j 2^{dj} \cdot 2^{-j\frac{d+1}{2}p},
\end{align} where the first inequality above follows from the fact that all norms on finite dimensional  normed spaces are comparable.  %\footnote{ When the denominator in (4.1) is bounded from below, then the integral over $\eta$ is bounded from above by the "surface area of unit sphere" which is a constant. I do not see why we have the  factors $2^{j(d+1)} \cdot 2^{-j} $  in the front of the inequality.}

\vskip.125in 

The conclusion follows after dividing both sides by $2^{dj}$ and taking $p$'th root. This recovers Theorem \ref{main} from Theorem \ref{maintool}. 

\vskip.25in 

\section{Sketch of proof of a slight strengthening of Theorem \ref{bourgain}} 
\label{bourgainth} 

\vskip.125in 

In this section we sketch the proof of the following result, which essentially follows from Bourgain's argument in \cite{Bou86}. 

\begin{definition} We say that $E \subset {\Bbb R}^d$ has positive upper Lebesgue density if 
$$ \overline{dens}_{{\mathcal L}}(E)=\limsup_{R \to \infty} \frac{|E \cap B_R|}{|B_R|}>0,$$ where $B_R$ denote the ball of radius $R$. 
\end{definition} 

\begin{definition} Suppose that $E \subset {\Bbb R}^d$  with  $\overline{dens}_{{\mathcal L}}(E)>0$. We say that $E'$ is a refinement of $E$ of relative density $r$ if $E' \subset E$ and 
$$\overline{dens}_{{\mathcal L}}(E') \ge r \cdot \overline{dens}_{{\mathcal L}}(E).$$
\end{definition} 

\begin{theorem} \label{bourgain+} Let $E \subset {\Bbb R}^d$ of positive upper Lebesgue density $\overline{dens}_{{\mathcal L}}(E)$. Then there exists a threshold $L_0>0$, a threshold $r \in (0,1)$, and $x_0 \in E$, such that for any refinement $E'$ of $E$ of relative density $r$, % \footnote{should not be the relative density be $r_0>r$ where $r$ is the threshold?}
 and any $L>L_0$, there exists $x' \in E'$ such that $\rho^{*}(x_0-x')=L$.
%\footnote{for example, if $\rho^{*}(x)= |x|$, and $E'$ is a small ball contained in a big ball $E$, then how can any big distances be realized?} 
 \end{theorem} 

\vskip.125in 

Bourgain's proof of the case $d=2$, $\rho^{*}(x)=|x|$, has two key ingredients. The first is the following maximal theorem which follows from Bourgain's proof of the circular maximal theorem. 

\begin{theorem} \label{2dmaximal} Let $\sigma$ be the arc-length measure on the circle $S^1$, and let $\sigma_s$ be defined by the relation $\widehat{\sigma}_s(\xi)=\widehat{\sigma}(s \xi)$. Let $P_t$ denote the Poisson semi-group kernel on ${\Bbb R}^2$, so $\widehat{P}_t(\xi)=e^{-t |\xi|}$. For $p>2$, there exist constants $C(p)<\infty$ and $\alpha(p)>0$, such that 
$$ {\left|\left| \sup_{s \ge t_0} |(f-(f*P_t)*\sigma_s)| \right|\right|}_p \leq C(p) {\left(\frac{t}{t_0} \right)}^{\alpha(p)} {||f||}_p.$$   

\end{theorem} 

The fact that the same conclusion is achieved if $S^1$ is replaced by any symmetric smooth convex curve with non-vanishing curvature follows from \cite{MSS93}. In higher dimensions, the following result follows from \cite{G81}. 

\begin{theorem} \label{allanmaximal} Let $\rho$ be the Minkowski functional on the convex symmetric body $K$ with a smooth boundary and everywhere non-vanishing curvature. Let $\rho^{*}$ be as above. Let $S_{\rho}=\{x \in {\Bbb R}^d: \rho^{*}(x)=1 \}$, and let $\sigma_{\rho}$ denotes the surface measure. Let 
$\sigma_{\rho,s}$ denote the dilated surface measure, as above. For $d>2$, there exist $C<\infty$, such that 
$$ {\left|\left| \sup_{s \ge t_0} |(f-(f*P_t)*\sigma_{\rho,s})| \right|\right|}_2 \leq C {\left(\frac{t}{t_0} \right)}^{\frac{d-2}{2}} {||f||}_2.$$ 
\end{theorem} 

\vskip.125in 

The second ingredient in Bourgain's proof is the observation that if the conclusion of Theorem \ref{bourgain} is false in the case $d=2, \rho^{*}(\xi)=|\xi|$, then there exists $U \subset {[0,1]}^2$, $|U|>\epsilon$, and a sequence of positive numbers 
$$ s_1>t_1>s_2>t_2>\dots>s_J>t_J,$$ $J$ arbitrarily large, satisfying 
$$ s_{j+1}<\frac{1}{2} t_j,$$ 

and $x \in U \cap [s_j, 1-s_j]$, $s_j \in (0 ,1/2)$, implies that $\sup_{s_j<t<s_{j+1}} |(1_{{[0,1]}^2}-f)*\sigma_t|=1$. The same reduction can be performed in dimensions $d \ge 3$ with the Euclidean metric replaced by $\rho^{*}$. The contradiction is then obtained using Theorem \ref{2dmaximal}. In dimensions three and higher, the contradiction is obtained in exactly the same way using Theorem \ref{allanmaximal}. 

Moreover, the same argument yields the slightly stronger conclusion, claimed in Theorem \ref{bourgain+}, that there exists a threshold $L_0>0$, a threshold $r \in (0,1)$, and $x_0 \in E$, such that for any refinement $E'$ of $E$ of relative density $r$, and any $L>L_0$, there exists $x' \in E'$ such that $\rho^{*}(x_0-x')=L$. This is because assuming that this statement is false leads to exactly the same statement as above and the contradiction is obtained in the same way. 

\vskip.25in 

\section{Proof of Theorem \ref{inftyforrealth}} 

\vskip.125in 

Suppose that ${\mathcal E}(A)$ is an orthogonal collection in $L^2(K)$ in the sense that 
$$ \int_K e^{2 \pi i x \cdot (a-a')} dx=0 \ \text{for all} \ a \not=a' \in A, $$ where $K$ is as in Theorem \ref{inftyforrealth}. 

The first listed author and Rudnev proved in \cite{IR03} that if $A$ is infinite, $A$ is contained in a line. Moreover, they showed that if $d \not=1 \mod 4$, then $A$ must be finite. 

The  proof relies on two estimates that follow from the orthogonality assumption for ${\mathcal E}(A)$. The first says that if ${\mathcal E}(A)$ is an orthogonal set in $L^2(K)$, then for all $a,a' \in A$, 
\begin{equation} \label{1pair} \rho^{*}(a-a')=\frac{k}{2}+\frac{d-1}{8}+O({|a-a'|}^{-1}) \ \text{for some integer} \ k \ \text{as} \  |a-a'| \to \infty.\end{equation}  

The second estimate says that if $a_0,a_1 \in A$, $a_0 \not=a_1$, then for all $a \in A$, 
\begin{equation} \label{2pair} \rho^{*}(a_0-a)-\rho^{*}(a_1-a)=\frac{k}{2}+O({|a|}^{-2}) \ \text{for some integer} \ k \ \text{as} \ |a| \to \infty. \end{equation} 

The first estimate follows directly from Lemma \ref{stationaryphase}. The second estimate uses more refined Bessel expansions (see Lemma 1.4 in \cite{IR03}). Using these estimates, the authors established their claim. 

\vskip.125in 

Under the hypotheses of Theorem \ref{inftyforrealth}, the same arguments yield the following form of (\ref{1pair}) and (\ref{2pair}): 
\begin{equation} \label{1pair_modified} \rho^{*}(a-a')=\frac{k}{2}+\frac{d-1}{8}+O \left( \frac{\phi(|a-a'|)}{{|a-a'|}^{\frac{d+1}{2}}} \right) 
\ \text{for some integer} \ k \ \text{as} \  |a-a'| \to \infty, \end{equation} where (\ref{previous}) holds, and (\ref{2pair}) becomes 
\begin{equation} \label{2pair_modified} \rho^{*}(a_0-a)-\rho^{*}(a_1-a)=\frac{k}{2}+O \left( \frac{\phi(|a|)}{{|a|}^{\frac{d+1}{2}}} \cdot {|a-a'|}^{-1} \right) 
\ \text{for some integer} \ k \ \text{as} \ |a| \to \infty. \end{equation} 

Observe that 
$$ O \left( \frac{\phi(|a-a'|)}{{|a-a'|}^{\frac{d+1}{2}}} \right) =o(1), \ \text{as} \ |a-a'| \to \infty,$$ by assumption. Similarly, 
$$ O \left( \frac{\phi(|a|)}{{|a|}^{\frac{d+1}{2}}} \cdot {|a-a'|}^{-1} \right)=o({|a-a|}^{-1}) \ \text{as} \ |a-a'| \to \infty.$$

The geometric observation used to establish the conclusion of the main result in \cite{IR03}, namely Lemma 1.5, does not require anything stronger than (\ref{1pair_modified}) and (\ref{2pair_modified}) above. These are all that are needed to establish Lemma 1.5 and Proposition 1.6 in that paper. This results in the following version of the Erd\H os Integer Distance Principle. 

\begin{theorem} \label{erdosnewth} (Approximate version of the Erd\H os Integer Distance Principle) Let $\rho^{*}$ be as above, corresponding to a symmetric convex body $K$ with a smooth boundary and non-vanishing Gaussian curvature. Let $A \subset {\Bbb R}^d$, $d \ge 2$, and suppose that for $a,a_0,a_1 \in A$, 

\begin{equation} \label{difference} \rho^{*}(a-a_0)-\rho^{*}(a-a_1)=o({|a|}^{-1}), \end{equation} and there exist constants $c_1,c_2$ such that for $a,a' \in A$, 

\begin{equation} \label{almostinteger} \rho^{*}(a-a')=c_1k+c_2+o(1). \end{equation} 

Then $A$ is a subset of a line. If $d \not=1 \mod 4$, then $A$ is finite. 

\end{theorem} 

\vskip.25in 

\section{Remarks and open problems} 

\vskip.125in 

The purpose of this section is to describe some open problems that may be accessible using the methods of this paper and related ideas. 

\vskip.125in 

\begin{itemize} 

\item[i)] The authors of this paper previously proved (\cite{IM2022}) that if $g(x)=\chi_{B_d}(x)$, where $B_d$ is the unit ball in dimensions two and higher, and $d \not=1 \mod 4$, then there does not exist $S \subset {\Bbb R}^d \times {\Bbb R}^d$, such that ${\{g(x-a)e^{2 \pi i x \cdot b} \}}_{(a,b) \in S}$ is an orthogonal basis for $L^2({\Bbb R}^d)$. The case when $d=1 \mod 4$ is work in progress. A reasonable question to ask is whether the same conclusion can be established under the assumption of approximate orthogonality. 

\vskip.125in 

\item[ii)] In view of \cite{IKT01}, it should be possible to show that the conclusion of Theorem \ref{main} still holds if we assume that $\partial K$ has one point, in the neighborhood of which it is smooth and has non-vanishing Gaussian curvature. 

\vskip.125in 

\item[iii)] It is reasonable to conjecture that if $K$ is a convex body in ${\Bbb R}^d$, $d \ge 2$, then there exists $S \subset {\Bbb R}^{2d}$ such that 
${\{\chi_K(x-a)e^{2 \pi i x \cdot b} \}}_{(a,b) \in S}$ is an orthogonal basis for $L^2(\Bbb R^d)$  if and only if $K$ tiles ${\Bbb R}^d$ by translation, i.e iff $K$ is a symmetric polytope of the type classified by Venkov \cite{V54} and McMullen \cite{M80}. It may be possible to approach this problem by combining the techniques of this paper and those of \cite{Lev-Mat} 
who recently proved that if $K \subset {\Bbb R}^d$, $d \ge 2$, is convex, then $L^2(K)$ has an orthogonal basis of exponentials if and only if $K$ tiles ${\Bbb R}^d$ by translation. The two dimensional case was previously established in \cite{IKT01}. 

\vskip.125in 

\item[iv)] The contrast between the proof of Theorem \ref{main} and Theorem \ref{inftyforrealth} raises the question of whether there is a strengthening of the approximate version of the Erd\H os Integer Distance Principle, proved in \cite{IR03}, where the point-wise (near) integer distance condition is replaced by an average version. For example, is it true that if $A$ is an infinite separated subset of ${\Bbb R}^d$, $d \ge 2$, such that there exist constants $c_1,c_2$ and $p \in [1,\infty)$ so that
$$ |a-a'|=c_1k+c_2+O({\mathcal D}(|a-a'|)),$$ with 
$$ \limsup_{j \to \infty} \frac{1}{2^j} \int_{2^j}^{2^{j+1}} {|{\mathcal D}(t)|}^p dt=0,$$ then $A$ is a subset of a line? 

If anything like this is true, it should lead to a considerable strengthening of Theorem \ref{main}. Heuristically, the idea (described in two dimensions for the sake of simplicity) is the following. Fix $a_0, a_1$, $a_0 \not=a_1$. Then $90$ percent of the points in $A$ are nearly $c_1k+c_2$ distance away from $a_0$, and the same is true with respect to $a_1$. It follows that $81$ percent of the points lie on finitely many hyperbolas with focal points at $a_0$ and $a_1$. Now choose $a_2$ that does not lie on a line that passes through $a_0$ and $a_1$ and play the same game with respect to $a_1, a_2$ and the points that are nearly shifted scaled integers with respect to $a_2$. The hyperbolas with focal points at $a_0,a_1$ intersect the hyperbolas with focal points at $a_1,a_2$ finitely many times (depending on the distance between those pairs of points) and shows that the resulting subset $A' \subset A$ (having at least $72.9$ percent of the points) is a subset of a line. We now replay the same game, using the same points $a_0,a_1,a_2$ with respect to the complement of $A'$. This approach can probably be used to show that $A=A_{good}+A_{bad}$, where $A_{good}$ is a subset of a line, and $A_{bad}$ is arbitrarily sparse. This may lead to an alternate proof of Theorem \ref{main}. Whether one can conclude that $A_{bad}$ is empty is an interesting question that should be investigated. 

\vskip.125in 

\end{itemize}


\vskip.25in 

\begin{thebibliography}{20}

\bibitem{B66} A. Beurling, {\it Local harmonic analysis with some applications to differential operators}, Some recent advances in the basic sciences, Academic press \textbf{1} (1966).

\bibitem{Bou86} J. Bourgain, {\it Szemeredi type theorem for sets of positive density in ${\Bbb R}^k$}, Israel J. Math. \textbf{54} (1986), no. 3, 307-316.

\bibitem{Bou86II} J. Bourgain, {\it Averages in the plane over convex curves and maximal operators}, J. Analyse Math. \textbf{47} (1986), 69--85.

%\bibitem{Fug01} B. Fuglede, {\it Orthogonal Exponentials on the Ball}, Expo. Math. \textbf{19}, (2001), 267-272. 
 
%\bibitem{Fug74} B. Fuglede, {\it Commuting self-adjoint partial differential operators and a group theoretic problem}, J. Funct. Anal. \textbf{16} (1974), 101-121.

\bibitem{FKW90} H. Furstenberg, Y. Katznelson, and B. Weiss, {\it Ergodic theory and configurations in sets of positive density} Mathematics of Ramsey theory, 184-198, Algorithms Combin., 5, Springer, Berlin, (1990).

\bibitem{GS58} I. Gelfand and G. Shilov, {\it Generalized Functions}, Vol. 1, Academic Press, (1958).

\bibitem{G81} A. Greenleaf, {\it Principal curvature and harmonic analysis}, Indiana Univ. Math. J. \textbf{30} (1981), no. 4, 519-537.

\bibitem{H62} C. Herz, {\it Fourier transforms related to convex sets}, Ann. of Math. (2) \textbf{75} (1962) 81-92. 

\bibitem{IK06} A. Iosevich and M. Kolountzakis, {\it A Weyl type formula for Fourier spectra and fraimes}, Proc. Amer. Math. Soc. \textbf{134} (2006), no. 11, 3267-3274.

%\bibitem{IK13} A. Iosevich and M. Kolountzakis, {\it Size of orthogonal sets of exponentials for the disk}, Rev. Mat. Iberoam. \textbf{29} (2013), no. 2, 739-747.

\bibitem{IM2020} A. Iosevich and A. Mayeli, {\it On complete and incomplete exponential systems},  to appear in  the Proceeding of AMS. 

\bibitem{IM2022} A. Iosevich and A. Mayeli, {\it Gabor orthogonal bases on convexity}, Discrete Anal. (2018), Paper No. \textbf{19}, 11 pp. 

\bibitem{IKP01} A. Iosevich, N.~H. Katz and S. Pedersen, {\it Fourier bases and a distance problem of Erd\"os}. \newblock {\em Math. Res. Lett.}, \textbf{6}, 251-255, (2001).

\bibitem{IKT01} A. Iosevich, N.~H. Katz and T. Tao, {\it Convex bodies with a point of curvature do not have Fourier bases}, Amer. J. Math., \textbf{123}, 115-120, (2001).

\bibitem{IR03} A. Iosevich and M. Rudnev, {\it A combinatorial approach to orthogonal exponentials}, Int. Math. Res. Not. (2003), no. 50, 2671-2685. 

\bibitem{IS96} A. Iosevich and E. Sawyer, {\it Oscillatory integrals and maximal averages over homogeneous surfaces}, Duke Math. J. \textbf{82} (1996), no. 1, 103-141.

\bibitem{IP00} A. Iosevich and S. Pedersen, {\it How large are the spectral gaps?} Pacific J. Math. \textbf{192} (2000), no. 2, 307-314.

%\bibitem{IR03} A. Iosevich and M. Rudnev, {\it A combinatorial approach to orthogonal exponentials}, Int. Math. Res. Not. (2003), no. 50, 2671-2685. 

 \bibitem{Lev-Mat} N. Lev, M. Matolcsi, {\it The Fuglede conjecture for convex domains is true in all dimensions}, Acta Mathematica, Volume 228 (2022), No. 2, Pages: 385 -- 420. 

\bibitem{M80} P. McMullen, {\it Convex bodies which tile space by translation}, Mathematika, \textbf{27} (1):113-121, (1980). 

\bibitem{MSS93} G. Mockenhaupt, A. Seeger, and C. D. Sogge, {\it Local smoothing of Fourier integral operators and Carleson-Sjolin estimates}, J. Amer. Math. Soc. \textbf{6} (1993), no. 1, 65-130.

\bibitem{St73} E. M. Stein, {\it  Maximal functions. I. Spherical means},  Proc. Nat. Acad. Sci. U.S.A. \textbf{73} (1976), no. 7, 2174-2175.

\bibitem{V54} B. A. Venkov, {\it On a class of Euclidean polyhedra}, Vestnik Leningrad. Univ. Ser. Mat. Fiz. Him., 9(2):11--31, (1954).

\end{thebibliography}
\end{document}